\documentclass[12pt]{article}
\usepackage{amsfonts}
\usepackage{mathrsfs}
\usepackage{amsmath,amssymb}
\baselineskip 4.5pt \lineskip 4.5pt \numberwithin{equation}{section}

\def\QED{\hfill$\Box$\par}

\def\HH{\mathcal {H}}

\def\LL{\mathcal {L}_{\lambda,\mu}}

\def\cl{\centerline}

\def\ni{\noindent}

\def\rar{\longrightarrow}

\def\vs{\vspace*}

\def\C{\mathbb{C}}

\def\Z{\mathbb{Z}}
\def\adddot{$\!\!\!${\bf.}\ \ }

\newtheorem{theo}{Theorem}[section]

\newtheorem{case}{Case}
\newtheorem{lemm}[theo]{Lemma}

\begin{document}
\baselineskip 18pt \cl{{\bf
 2-Cocycles of Original Deformative
Schr\"{o}dinger-Virasoro Algebras}\footnote {Supported by NSF grants
10471091, 10671027 of China, ``One Hundred Talents Program'' from
University of Science
and Technology of China.\\[2pt] \indent Corresponding E-mail:
sd\_junbo@163.com}} \vs{6pt}

\cl{Junbo Li$^{*,\dag)}$, Yucai Su$^{\ddag)}$}

\cl{\small $^{*)}$Department of Mathematics, Shanghai Jiao Tong
University, Shanghai 200240, China}

\cl{\small $^{\dag)}$ Department of Mathematics, Changshu Institute
of Technology, Changshu 215500, China}

\cl{\small $^{\ddag)}$ Department of Mathematics, University of
Science and Technology of China, Hefei 230026, China}

\cl{\small E-mail: sd\_junbo@163.com, ycsu@ustc.edu.cn}\vs{6pt}

\noindent{\small{\bf Abstract.} Both original and twisted
Schr\"{o}dinger-Virasoro algebras also their deformations were
introduced and investigated in a series of papers by Henkel, Roger
and Unterberger. In the present paper we aim to determine the
2-cocycles of original deformative Schr\"{o}dinger-Virasoro
algebras.

} \vs{12pt}

\cl{\bf\S1. \
Introduction}\setcounter{section}{1}\setcounter{theo}{0}\setcounter{equation}{0}

Both original and twisted Schr\"{o}dinger-Virasoro algebras also
their deformations were introduced in \cite{H1,H2,HU,RU}, in the
context of non-equilibrium statistical physics, closely related to
both Schr\"{o}dinger Lie algebras and the Virasoro Lie algebra,
which are known to be important in many areas of mathematics and
physics (e.g., statistical physics). Their vertex representations
were constructed in \cite{U}. Later the derivation algebra and
automorphism group of the twisted sector were determined in
\cite{LS1}. Almost at the same time, the derivations, central
extensions and automorphism group of the extended sector were
investigated in \cite{GJP}. Furthermore, irreducible modules with
finite-dimensional weight spaces and indecomposable modules over
both original and twisted sectors were considered in \cite{LS2}, and
the second cohomology group of a class of twisted deformative
sectors was determined in \cite{LS3} by the authors.

The infinite-dimensional Lie algebras $\LL\,(\lambda,\mu\in\C)$
considered in this paper called {\it original deformative
Schr\"{o}dinger-Virasoro Lie algebras} (\,see \cite{RU}), possess
the same $\C$-basis
$$\{L_n,\,M_n,\,Y_p\,|\,n\in \Z,\,p\in\frac{1}{2}+\Z\}$$
with the following Lie brackets:
\begin{eqnarray}
&[L_n,L_m]\!\!\!&=(m-n)L_{m+n},\label{LB1}
\\[4pt]
&[L_n,Y_m]\!\!\!&=(m-\frac{(\lambda+1)n}{2}+\mu)Y_{m+n},\ \ \,\,
[Y_n,Y_{m}]=(m-n)M_{m+n},\label{LB2}\\[4pt]
&[L_n,M_m]\!\!\!&=(m-\lambda n+2\mu)M_{m+n},\ \ \ \,\ \ \ \ \,\,
[Y_n,M_m]=[M_n,M_m]=0.\label{LB3}
\end{eqnarray}

The purpose of this paper is to determine the 2-cocycles of the
original deformative Schr\"{o}dinger-Virasoro algebras
$\LL\,(\lambda,\mu\in\C)$ defined above. The 2-cocycles on Lie
algebras play important roles in the central extensions of Lie
algebras, which can be used to construct many infinite-dimensional
Lie algebras, and further to describe the structures and some of the
representations of these type Lie algebras. As the cohomology groups
are closely related to the structure of Lie algebras, the
computation of cohomology groups seems to be important and
interesting as well. Maybe due to the reasons stated above, there
appeared a number of papers on 2-cocycles and cohomology groups of
infinite-dimensional Lie algebras and conformal algebras (\,see
\cite{BKV,GJP}, \cite{L,LW} and \cite{SZR}--\cite{SZK}). Now let's
formulate our main results below.

Recall that a {\it 2-cocycle} on some $\LL$ is a $\C$-bilinear
function $\psi:\LL\times \LL\rar \C$ satisfying the following
conditions:
\begin{eqnarray}
&&\psi(v_1,v_2)=-\psi(v_2,v_1)\mbox{\ \ (\,skew-symmetry)},\nonumber\\
&&\psi([v_1,v_2],v_3)+\psi([v_2,v_3],v_1)+\psi([v_3,v_1],v_2) =0
\mbox{\ \ (Jacobian identity)}\label{2c-Ji},
\end{eqnarray}
for $v_1,v_2,v_3\in \LL$. Denote the vector space of 2-cocycles on
$\LL$ by $\mathcal {C}^2(\LL,\C)$. For any $\C$-linear function
$f:\LL\rar\C$, define a 2-cocycle $\psi_f$ as follows
\begin{eqnarray}
\psi_f(v_1,v_2)=f([v_1,v_2]),\ \ \,\forall\;v_1,v_2\in
\LL,\label{cobco}
\end{eqnarray}
which is usually called a {\it 2-coboundary} or a {\it trivial
2-cocycle} on $\LL$. Denote the vector space of 2-coboundaries on
$\LL$ by $\mathcal {B}^2(\LL,\C)$. A 2-cocycle $\varphi$ is said to
be {\it equivalent to} a 2-cocycle $\psi$ if $\varphi-\psi$ is
trivial. For a 2-cocycle $\psi$, we denote the equivalent class of
$\psi$ by $[\psi]$. The quotient space $\mathcal
{H}^2(\LL,\C)\!=\!\mathcal {C}^2(\LL,\C)/\mathcal {B}^2(\LL,\C)$ is
called the {\it second cohomology group} of $\LL$.

Usually one calls a 2-cocycle $\xi$ on $\LL$ the {\it Virasoro
cocycle}, denoted by $\xi_{Vir}$, if
\begin{eqnarray}\label{defvirco}
&&\xi(L_n,L_{m})=\frac{n^3-n}{12}\delta_{m,-n}\,,\ \ {\mbox {while
other components vanishing}}.
\end{eqnarray}

The main results of the paper can be formulated as follows.
\begin{theo}\adddot\label{mainth}
(i)\ If $\mu\notin\{\frac{1}{2}\Z\}$, then for any $\lambda\in\C$,
$\HH^2(\mathcal {L}{_{\lambda,\mu}},\C)\cong\C$ is generated
by the Virasoro cocycle.\\
(ii)\ For $\mu\in\frac{1}{2}+\Z$ and $\lambda\ne-3,-1,1$,
$\HH^2(\mathcal {L}{_{\lambda,\mu}},\C)\cong\C$ is generated
by the Virasoro cocycle.\\
(iii)\ For $\mu\in\frac{1}{2}+\Z$ and $\lambda=-3$, $\HH^2(\mathcal
{L}_{-3,\mu},\C)\cong\C^2$ is generated by the Virasoro cocycle and
an independent cocycle of the form $c(L_n,Y_m)=\delta_{n,-m-\mu}$
(all other
components vanishing).\\
(iv)\ For $\mu\in\frac{1}{2}+\Z$ and $\lambda=-1$, $\HH^2(\mathcal
{L}_{-1,\mu},\C)\cong\C^3$ is generated by the Virasoro cocycle and
other two independent cocycles $c_1$ and $c_2$ defined by (all other
components vanishing)
\begin{eqnarray*}
c_1(M_m,Y_{n})=\delta_{n,-m-3\mu},\ \ c_2(L_{-m},Y_{n})
=\frac{m(m+1)}{2}\delta_{n,m-\mu}.
\end{eqnarray*}
(v)\ For $\mu\in\frac{1}{2}+\Z$ and $\lambda=1$, $\HH^2(\mathcal
{L}_{1,\mu},\C)\cong\C^3$ is generated by the Virasoro cocycle and
other two independent cocycles $c_1$ and $c_2$ defined by (all other
components vanishing)
\begin{eqnarray*}
&&c_1(L_{-m},Y_{n})=m(m^2-1)\delta_{n,m-\mu},\\
&&c_2(L_{-m},M_{n-2\mu})=c_2(Y_{-m-\mu},Y_{n-\mu})=m(m^2-1)\delta_{m,n}.
\end{eqnarray*}
(vi)\ For $\mu\in\Z$, if $\lambda\ne-1$, $\HH^2(\mathcal
{L}{_{\lambda,\mu}},\C)\cong\C$ is generated by the Virasoro
cocycle; while if $\lambda=-1$, $\HH^2(\mathcal
{L}_{-1,\mu},\C)\cong\C^2$ is generated by the Virasoro cocycle and
an independent cocycle of the form
$c(Y_{p},Y_{q})=-\delta_{q,-p-2\mu}(p+\mu)$ (all other components
vanishing).
\end{theo}

Throughout the paper, we denote by $\Z^*$ the set of all nonzero
integers, $\C^*$ the set of all nonzero complex numbers and
$\C^*\!\!\setminus\!\Z^*=\{x\,|\,x\in\C^*,\,x\notin\Z^*\}$.

\vs{18pt} \cl{\bf\S2. \ Proof the main
results}\setcounter{section}{2}\setcounter{theo}{0}\setcounter{equation}{0}

Let $\psi$ be any 2-cocycle. Define a $\C$-linear function
$f:\LL\rightarrow\C$ as follows
\begin{eqnarray}
f(L_n)\!\!\!&=&\!\!\!\left\{\begin{array}{ll}
\frac{1}{n}\psi(L_0,L_n) &{\rm if}\ \,n\neq0,\,\,\forall\,\,\mu\in\C^*,\vs{6pt}\\
\frac{1}{2}\psi(L_{-1},L_{1}) &{\rm if}\
\,n=0,\,\,\forall\,\,\mu\in\C^*,
\end{array}\right.\label{delfL}\\
f(M_n)\!\!\!&=&\!\!\!\left\{\begin{array}{ll}
\frac{1}{n+2\mu}\psi(L_0,M_n) &{\rm if}\
\,n\neq-2\mu,\,\mu\in\frac{1}{2}\Z,\ {\rm or}\
\mu\notin\frac{1}{2}\Z,
\vs{6pt}\\
\frac{-1}{\lambda+1}\psi(L_1,M_{-2\mu-1}) &{\rm if}\
\,n=-2\mu,\,\lambda\neq-1, \ {\rm
and}\ \mu\in\frac{1}{2}\Z,\end{array}\right.\label{delfM}\\
f(Y_p)\!\!\!&=&\!\!\!\left\{\begin{array}{ll}
\frac{1}{p+\mu}\psi(L_0,Y_p) &{\rm if}\ \,p\neq-\mu,\
\mu\in\frac{1}{2}+\Z,\
{\rm or}\ \mu\notin\frac{1}{2}+\Z,\vs{6pt}\\
\frac{-2}{\lambda+3}\psi(L_1,Y_{-\mu-1}) &{\rm if}\ \,p=-\mu,
\lambda\neq-3\ {\rm and}\
\mu\in\frac{1}{2}+\Z.\end{array}\right.\label{delfY}
\end{eqnarray}
Let $\varphi=\psi-\psi_f-\xi_{Vir}$ where $\psi_f$ and $\xi_{Vir}$
are respectively defined in (\ref{cobco}) and (\ref{defvirco}), then
\begin{eqnarray}\label{recl1}
\varphi(L_m,L_n)=0,\ \ \forall\,\,m,\,n\in\Z.
\end{eqnarray}

The proof of Theorem \ref{mainth} is included in the following three
technical lemmas.
\begin{case}
\adddot\label{cs1} $\mu\notin\frac{1}{2}\Z$.
\end{case}
\begin{lemm}\adddot\label{lemm1}
$\varphi=0$.
\end{lemm}
{\it Proof.} According to (\ref{delfM}) and (\ref{delfY}), one has
\begin{eqnarray}\label{recl18998}
\varphi(L_0,Y_p)=\varphi(L_0,M_n)=0,\ \
\forall\,\,n\in\Z,\,p\in\frac{1}{2}+\Z.
\end{eqnarray}
For any $p,q\in\frac{1}{2}+\Z$, using the Jacobian identity on the
triple $(L_0,\,Y_{p},\,Y_q)$, together with (\ref{recl1}), we obtain
\begin{eqnarray*}
(p+q+2\mu)\varphi(Y_p,Y_q)=0,
\end{eqnarray*}
which together with our assumption $\mu\notin\frac{1}{2}\Z$ forces
\begin{eqnarray}\label{corI-4}
\varphi(Y_p,Y_q)=0.
\end{eqnarray}
For any $m,n\in\Z,\,p\in\frac{1}{2}+\Z$, using the Jacobian identity
on the four triples $(L_{m},\,Y_p,\,L_0)$, $(L_{m},\,M_n,\,L_0)$,
$(Y_{p},\,M_n,\,L_0)$ and $(M_{m},\,M_n,\,L_0)$ in (\ref{2c-Ji})
respectively, one has
\begin{eqnarray}
&&(\,m+p+\mu\,)\,\varphi(\,L_m,Y_p\,)\,=0,\label{corIad-1}\\
&&(m+n+2\mu)\,\varphi(L_m,M_n)=0,\label{corIad-2}\\
&&(p+n+3\mu)\,\varphi(Y_p,M_n)\,=0,\label{corI-5}\\
&&(m+n+4\mu)\varphi(M_m,M_n)=0.\label{corI-6}
\end{eqnarray}
Then our assumption $\mu\notin\frac{1}{2}\Z$, together with
(\ref{corIad-1}) and (\ref{corIad-2}) gives
\begin{eqnarray}
\varphi(L_m,Y_p)=\varphi(Y_p,M_n)=\varphi(L_m,M_n)=0.\label{corIad-3}
\end{eqnarray}
Similar to the proof of Subcase 1.3 given in \cite{LS1}, one also
has
\begin{eqnarray}
&&\varphi(M_m,M_n)=0,\ \ \forall\,\,m,n\in\Z.\label{corI-I-3}
\end{eqnarray}
Then this lemma follows.\QED

\ni This lemma in particular proves Theorem \ref{mainth} (i).
\begin{case}
\adddot\label{cs2} $\mu\in\frac{1}{2}+\Z$.
\end{case}
\begin{lemm}\adddot\label{lemm2}
For any $\lambda\in\C$, $\varphi=0$ unless the following subcases:

i)\ \,if $\lambda=-1$, then $\varphi(M_m,Y_{-m-3\mu})=
\varphi(M_0,Y_{-3\mu})$, $\varphi(L_{-m},Y_{m-\mu})
=\frac{m(m+1)}{2}\varphi(L_{-1},Y_{1-\mu})$;

ii)\ \,if $\lambda=-3$, then
$\varphi(L_{-m},Y_{m-\mu})=\varphi(L_{0},Y_{-\mu})$;

iii)\ \,if $\lambda=1$, then
$\varphi(L_{-m},Y_{m-\mu})=m(m^2-1)c_1$,
$\varphi(Y_{-m-\mu},Y_{m-\mu})=m(m^2-1)c'_1$,
$\varphi(L_{-m},M_{m-2\mu})=m(m^2-1)c'_1$;

\ni for any $m\in\Z$, and some constants
$c_1,\,c'_1\in\C$.
\end{lemm}
{\it Proof.} According to (\ref{delfM})
and (\ref{delfY}), one has
\begin{eqnarray}
&&\varphi(L_0,Y_p)=0\ \ {\rm if}\ \,p\neq-\mu,\ \
\varphi(L_1,Y_{-\mu-1})=0\ \ {\rm if}\ \,\lambda\neq-3,\label{recII2}\\
&&\varphi(L_0,M_n)=0\ \ {\rm if}\ \,n\neq-2\mu,\ \
\varphi(L_1,M_{-2\mu-1})=0\ \ {\rm if}\
\,\lambda\neq-1.\label{recII1}
\end{eqnarray}
For any $m,n\in\Z,\,p,q\in\frac{1}{2}+\Z$, using the Jacobian
identity on the triples $(L_0,\,Y_{p},\,M_m)$,
$(L_0,\,M_{m},\,M_n)$, $(L_0,\,Y_{p},\,Y_q)$, $(L_0,\,L_{m},\,M_n)$
and $(L_0,\,L_{m},\,Y_p)$ in (\ref{2c-Ji}) respectively, together
with (\ref{recII2}) and (\ref{recII1}), one has
\begin{eqnarray}
&&(p+m+3\mu)\varphi(Y_p,M_m)=0,\label{corIIYM}\\
&&(m+n+4\mu)\varphi(M_m,M_n)=0,\label{corIIMM}\\
&&(p-q)\varphi(L_0,M_{p+q})+(p+q+2\mu)\varphi(Y_p,Y_q)=0,\label{corIIYY}\\
&&(m\lambda-2\mu-n)\varphi(L_0,M_{m+n})+(m+n+2\mu)\varphi(L_m,M_n)=0,\label{corIILM}\\
&&\big(m(1+\lambda)-2(p+\mu)\big)\varphi(L_0,Y_{m+p})
+2(m+p+\mu)\varphi(L_m,Y_p)=0.\label{corIILY}
\end{eqnarray}
For any $m\in\Z,\,p\in\frac{1}{2}+\Z,\,p\ne2\mu-\frac{m}{2}$, using
the Jacobian identity on the triple $(Y_{p-4\mu},\,Y_{-p-m},\,M_m)$,
one has
\begin{eqnarray*}
&&(m+2p-4\mu)\varphi(M_{-m-4\mu},M_m)=0,
\end{eqnarray*}
which together with (\ref{corIIMM}), gives
\begin{eqnarray}
&&\varphi(M_m,M_n)=0,\ \ \forall\,\,m,n\in\Z.\label{corIIMMla}
\end{eqnarray}
Using (\ref{recII1}), the identities (\ref{corIIYY}) and
(\ref{corIILM}) can be rewritten as
\begin{eqnarray}
\varphi(Y_p,Y_q)=\varphi(L_0,M_{-2\mu})=\varphi(L_m,M_n)=0\ \ {\rm
if}\ \,(p+q+2\mu)(m+n+2\mu)\ne0.\label{corIIYYLM-1}
\end{eqnarray}
The identities (\ref{recII2}) and (\ref{corIILY}) together give
\begin{eqnarray}
&&\varphi(L_m,Y_p)=0\ \ {\rm if}\ \,m+p\neq-\mu,\ \ \
\varphi(\,L_0,Y_{-\mu}\,)=0\ \ \ {\rm if}\
\,\lambda\neq-3.\label{corIILY-1}
\end{eqnarray}
From (\ref{recl1}), (\ref{corIIYM}) and
(\ref{corIIMMla})--(\ref{corIILY-1}), the left components we have to
present in this case are listed in the following
(\,$\forall\,\,n\in\Z,\,p\in\frac{1}{2}+\Z$):
\begin{eqnarray*}
\varphi(Y_n,L_{-\mu-n}),\ \,\varphi(L_{-2\mu-n},M_n),\
\,\varphi(Y_{-2\mu-p},Y_p)\ \,{\rm and}\ \,\varphi(M_n,Y_{-n-3\mu}).
\end{eqnarray*}
\ni{\bf Step 1.}\ \,The computation of
$\varphi(M_n,Y_{-n-3\mu}),\,\,\forall\,\,n\in\Z$.

For any $m,n\in\Z$, using the Jacobian identity on the triple
$(L_{n-m},Y_{-n-3\mu},M_m)$, one has
\begin{eqnarray}
&&2\big(m(1+\lambda)+2\mu-n\lambda\big)\varphi(M_n,Y_{-n-3\mu}\nonumber)\\
&&=\big(4\mu+(3+\lambda)n-(\lambda+1)m\big)\varphi(M_m,Y_{-m-3\mu}).\label{Lem2.1-1}
\end{eqnarray}
If $\lambda=0$, taking $m=0$ in (\ref{Lem2.1-1}), one has
\begin{eqnarray}
\varphi(M_n,Y_{-n-3\mu})
=(1+\frac{3n}{4\mu})\varphi(M_0,Y_{-3\mu}).\label{Lem2.1-1.1}
\end{eqnarray}
Taking (\ref{Lem2.1-1.1}) back to (\ref{Lem2.1-1}), one has
\begin{eqnarray*}
(m-n)\big(m(1+\lambda)+2\lambda(n+2\mu)\big)\varphi(M_0,Y_{-3\mu})=0,
\end{eqnarray*}
which forces $\varphi(M_0,Y_{-3\mu})=0$ and further
$\varphi(M_n,Y_{-n-3\mu})=0$, $\forall\,\,n\in\Z$
by (\ref{Lem2.1-1.1}).\\
If $\lambda\ne0,\,\frac{2\mu}{\lambda}\notin\Z$, taking $m=0$ in
(\ref{Lem2.1-1}), one has
\begin{eqnarray*}
\varphi(M_n,Y_{-n-3\mu})
=\frac{3n+4\mu+n\lambda}{4\mu-2n\lambda}\varphi(M_0,Y_{-3\mu}).
\end{eqnarray*}
Similarly, one also can prove
$\varphi(M_n,Y_{-n-3\mu})=\varphi(M_0,Y_{-3\mu})=0$.\\
If $\lambda\ne0,\,\frac{2\mu}{\lambda}\in\Z$, taking $m=0$ in
(\ref{Lem2.1-1}), one has
\begin{eqnarray}
\varphi(M_n,Y_{-n-3\mu})
=\frac{3n+4\mu+n\lambda}{4\mu-2n\lambda}\varphi(M_0,Y_{-3\mu})\ \ \
{\rm if}\ \,n\neq\frac{2\mu}{\lambda}.\label{Lem2.1-1.2}
\end{eqnarray}
Taking $(m,n,\lambda)=(\frac{2\mu}{\lambda},1,1)$ and
$(m,n)=(\frac{2\mu}{\lambda},0)$ in (\ref{Lem2.1-1}) respectively,
and using (\ref{Lem2.1-1.2}), one has
\begin{eqnarray}\label{Lemmaa1}
\varphi(M_{\frac{2\mu}{\lambda}},Y_{-\frac{2\mu}{\lambda}-3\mu})=\left\{\begin{array}{lll}
\frac{(6\mu-1)(1+\mu)}{2\mu-1}\varphi(M_0,Y_{-3\mu}) &{\rm if}\ \,\lambda=1,\,\mu\ne\frac{1}{2},\vs{6pt}\\
\varphi(M_1,Y_{-\frac{5}{2}}) &{\rm if}\ \,\lambda=1,\,\mu=\frac{1}{2},\vs{6pt}\\
\frac{2(1+2\lambda)}{\lambda-1}\varphi(M_0,Y_{-3\mu}) &{\rm if}\
\,\lambda\ne0,1,\, \frac{2\mu}{\lambda}\in\Z.
\end{array}\right.\label{delfL}
\end{eqnarray}
For the special case $\lambda=1,\,\mu=\frac{1}{2}$, by taking $n=2$
in (\ref{Lem2.1-1.2}) and $m=2,\,n=1$ in (\ref{Lem2.1-1}), we obtain
\begin{eqnarray*}
\varphi(M_2,Y_{-\frac{7}{2}})=4\varphi(M_1,Y_{-\frac{5}{2}})=-5\varphi(M_0,Y_{-\frac{3}{2}}),
\end{eqnarray*}
which gives
\begin{eqnarray}\label{lema6.12}
\varphi(M_1,Y_{-\frac{5}{2}})=-\frac{5}{4}\varphi(M_0,Y_{-\frac{3}{2}}).
\end{eqnarray}
Taking $(m,n)=(0,\frac{2\mu}{\lambda})$ in (\ref{Lem2.1-1}), one has
$\frac{3\mu(\lambda+1)}{\lambda}\varphi(M_0,Y_{-3\mu})=0$, which
infers
\begin{eqnarray}\label{lema6.16}
\varphi(M_0,Y_{-3\mu})=0\ \ \ {\rm if}\ \,\lambda\ne0,-1.
\end{eqnarray}
Then adding up the results obtained after (\ref{lema6.12}), one has
\begin{eqnarray}\label{Lemma2-1main}
\varphi(M_n,Y_{-n-3\mu})=\left\{\begin{array}{cl}
0&{\rm if}\ \,-1\ne\lambda\in\C,\,\mu\in\frac{1}{2}+\Z,\vs{6pt}\\
\varphi(M_0,Y_{-3\mu})&{\rm if}\
\,\lambda=-1,\,\mu\in\frac{1}{2}+\Z.
\end{array}\right.
\end{eqnarray}

\ni{\bf Step 2.}\ \,The computation of
$\varphi(Y_n,L_{-\mu-n}),\,\,\forall\,\,n\in\Z$.

Applying the Jacobian identity on the triple
$(L_{-m},L_n,Y_{m-n-\mu}),\,\,\forall\,\,m,n\in\Z$, one has
\begin{eqnarray}\label{L1S21}
&&\big(2m-n(\lambda+3)\big)\varphi(L_{-m},Y_{m-\mu})\nonumber\\
&&=\big(m(\lambda+3)-2n\big)\varphi(L_n,Y_{-n-\mu})
+2(m+n)\varphi(L_{n-m},Y_{m-n-\mu}).
\end{eqnarray}
Replacing $n$ by $-n$ and $m$ by $m+n$ in (\ref{L1S21})
respectively, one has
\begin{eqnarray}
&&\big(2m+n(\lambda+3)\big)\varphi(L_{-m},Y_{m-\mu})\nonumber\\
&&=\big(m(\lambda+3)+2n\big)\varphi(L_{-n},Y_{n-\mu})
+2(m-n)\varphi(L_{-n-m},Y_{m+n-\mu}),\label{L1S22}\\
&&\big(2(m+n)-n(\lambda+3)\big)\varphi(L_{-m-n},Y_{m+n-\mu})\nonumber\\
&&=\big((m+n)(\lambda+3)-2n\big)\varphi(L_n,Y_{-n-\mu})
+2(m+2n)\varphi(L_{-m},Y_{m-\mu}).\label{L1S23}
\end{eqnarray}
If $\lambda=-3$, then (\ref{L1S21}) can be rewritten as (\,by taking
$n=1$)
\begin{eqnarray}\label{L1S21}
&&(m+1)\varphi(L_{1-m},Y_{m-1-\mu})-m\varphi(L_{-m},Y_{m-\mu})=\varphi(L_1,Y_{-1-\mu}).
\end{eqnarray}
Using induction on $m$ in (\ref{L1S21}), we obtain
\begin{eqnarray}\label{L1S220}
\varphi(L_{-m},Y_{m-\mu})=\varphi(L_{0},Y_{-\mu}),\ \
\forall\,\,m\in\Z.
\end{eqnarray}
If $\lambda\ne-3$, taking $n=-1$ in both (\ref{L1S22}) and
(\ref{L1S23}), together with (\ref{recII2}), we have
\begin{eqnarray}
&&\big(2m-(\lambda+3)\big)\varphi(L_{-m},Y_{m-\mu})
=2(m+1)\varphi(L_{1-m},Y_{m-1-\mu}),\label{L1S24}\\
&&\big(2(m-1)+(\lambda+3)\big)\varphi(L_{1-m},Y_{m-1-\mu})\nonumber\\
&&=\big((m-1)(\lambda+3)+2\big)\varphi(L_{-1},Y_{1-\mu})
+2(m-2)\varphi(L_{-m},Y_{m-\mu}).\label{L1S25}
\end{eqnarray}
If $(\lambda+5)(\lambda-1)\ne0$, then identities (\ref{L1S24}) and
(\ref{L1S25}) together give
\begin{eqnarray}
&&\varphi(L_{-m},Y_{m-\mu})
=\frac{2(m+1)\big((\lambda+1)-m(\lambda+3)\big)}{{\lambda}^2+4\lambda-5}\varphi(L_{-1},Y_{1-\mu}),\label{L1S26}
\end{eqnarray}
Taking $\varphi(L_{-m},Y_{m-\mu})$ and
$\varphi(L_{1-m},Y_{m-1-\mu})$ obtained from (\ref{L1S26}) back to
(\ref{L1S24}), one has
\begin{eqnarray}\label{eqbu1}
\frac{m(2m-\lambda-3)({\lambda}^2+4\lambda+3)}{{\lambda}^2+4\lambda-5}\varphi(L_{-1},Y_{1-\mu})=0,\
\ \forall\,\,m\in\Z,
\end{eqnarray}
which forces (\,since the index $m$ can be shifted and our
assumption $\lambda\ne-3$)
\begin{eqnarray*}
({\lambda}^2+4\lambda+3)\varphi(L_{-1},Y_{1-\mu})=0.
\end{eqnarray*}
In another word,
\begin{eqnarray}\label{cali1}
&&\!\!\!\!\!\!\!\!\!\!\!\!\!\!\!\!{\mbox{the system consisted of
linear equations (\ref{L1S24}) and
(\ref{L1S25}) has nonzero}}\nonumber\\
&&\!\!\!\!\!\!\!\!\!\!\!\!\!\!\!\!{\mbox{solutions if and only if
$\lambda=-1$ under our assumption
$(\lambda+5)(\lambda-1)(\lambda+3)\ne0$.}}
\end{eqnarray}
If $\lambda=-1$, then (\ref{L1S26}) can be rewritten as
\begin{eqnarray}
&&\varphi(L_{-m},Y_{m-\mu})
=\frac{m(m+1)}{2}\varphi(L_{-1},Y_{1-\mu}),\ \
\forall\,\,m\in\Z.\label{L1S26-1}
\end{eqnarray}
If $\lambda=1$, then (\ref{L1S24}) becomes
\begin{eqnarray}
&&(m-2)\varphi(L_{-m},Y_{m-\mu})
=(m+1)\varphi(L_{1-m},Y_{m-1-\mu}),\label{L1S27}
\end{eqnarray}
which further gives \big(\,by taking $m=2$ in (\ref{L1S27})\big)
\begin{eqnarray}\label{L1S28}
\varphi(L_{-1},Y_{1-\mu})=0.
\end{eqnarray}
Then using induction on $m$ in (\ref{L1S27}), one can deduce
\begin{eqnarray}\label{L1S29}
\varphi(L_{-m},Y_{m-\mu})=\left\{\begin{array}{ll}
m(m^2-1)c_1&{\rm if}\ \,m\geq-2\mu-1,\vs{6pt}\\
m(m^2-1)c_2 &{\rm if}\ \,m<-2\mu-1,
\end{array}\right.
\end{eqnarray}
for some constants $c_1,\,c_2\in\C$. One thing left to be done is to
investigate the relations between the constants $c_1$ and $c_2$. If
$\lambda=1$, then (\ref{L1S25}) becomes
\begin{eqnarray*}
&&(m+2n)\varphi(L_{-m},Y_{m-\mu})\nonumber\\
&&=(2m+n)\varphi(L_{-n},Y_{n-\mu})
+(m-n)\varphi(L_{-n-m},Y_{m+n-\mu}),\\
&&(m-n)\varphi(L_{-m-n},Y_{m+n-\mu})\nonumber\\
&&=(2m+n)\varphi(L_n,Y_{-n-\mu}) +(m+2n)\varphi(L_{-m},Y_{m-\mu}).
\end{eqnarray*}
which together with each other force
\begin{eqnarray}
\big(2m+n\big)\label{L1S210}
\big(\varphi(L_{-n},Y_{n-\mu})+\varphi(L_{n},Y_{-n-\mu})\big)=0,
\end{eqnarray}
Then using (\ref{L1S29}) and (\ref{L1S210}), we obtain $c_1=c_2$,
which together with (\ref{L1S29}) gives
\begin{eqnarray}\label{L1S211}
\varphi(L_{-m},Y_{m-\mu})=m(m^2-1)c_1,\ \ \forall\,\,m\in\Z.
\end{eqnarray}
If $\lambda=-5$, then (\ref{L1S22}) and (\ref{L1S23}) convert to the
following form:
\begin{eqnarray}
&&(m-n)\varphi(L_{-m},Y_{m-\mu})\nonumber\\
&&=(n-m)\varphi(L_{-n},Y_{n-\mu})
+(m-n)\varphi(L_{-n-m},Y_{m+n-\mu}),\label{L1S2-51}\\
&&(m+2n)\varphi(L_{-m-n},Y_{m+n-\mu})\nonumber\\
&&=-(m+2n)\varphi(L_n,Y_{-n-\mu})
+(m+2n)\varphi(L_{-m},Y_{m-\mu}).\label{L1S2-52}
\end{eqnarray}
Furthermore, taking $n=-1$ in both (\ref{L1S2-51}) and
(\ref{L1S2-52}), and using (\ref{recII2}), one has
\begin{eqnarray*}
&&(m+1)\big(\varphi(L_{-m},Y_{m-\mu})-\varphi(L_{1-m},Y_{m-1-\mu})\big)=0,\\
&&(m-2)\big(\varphi(L_{-m},Y_{m-\mu})-\varphi(L_{1-m},Y_{m-1-\mu})\big)
=(m-2)\varphi(L_{-1},Y_{1-\mu}).
\end{eqnarray*}
from which and using (\ref{recII2}) again, can we deduce the
following relation:
\begin{eqnarray}\label{relalast0}
\varphi(L_{-m},Y_{m-\mu})=0,\ \ \forall\,\,m\in\Z.
\end{eqnarray}
If $\lambda\notin\{-5,\,-3,\,-1,\,1\}$, the identities (\ref{L1S26})
and (\ref{eqbu1}) together also force (\ref{relalast0}) to hold.

\ni{\bf Step 3.}\ \,The computation of
$\varphi(L_{-2\mu-n},M_n),\,\,\forall\,\,n\in\Z$.

Using the Jacobian identity on
$(L_{-2\mu-m},\,L_{n},\,M_{m-n}),\,\,\forall\,\,m,n\in\Z$, we obtain
\begin{eqnarray}
&&\!\!\!\!\!\!\!\!\!\big(m+2\mu-n(1+\lambda)\big)\varphi(L_{-m-2\mu},M_{m})\nonumber\\
&&\!\!\!\!\!\!\!\!\!=\big((m+2\mu)(1+\lambda)-n\big)\varphi(L_{n},M_{-n-2\mu})+(m+2\mu+n)
\varphi(L_{n-m-2\mu},M_{m-n}).\label{corII.1ad-n}
\end{eqnarray}
Replacing $n$ by $-n$ and $m$ by $m+n$ in (\ref{corII.1ad-n}), one
respectively gets
\begin{eqnarray}
&&\!\!\!\!\!\!\!\!\!\!\!\!\!\!\!\!\!\!\!
\big(m+2\mu+n(1+\lambda)\big)\varphi(L_{-m-2\mu},M_{m})\nonumber\\
&&\!\!\!\!\!\!\!\!\!\!\!\!\!\!\!\!\!\!\!
=(m+2\mu-n)\varphi(L_{-n-m-2\mu},M_{m+n})+
\big((m+2\mu)(1+\lambda)+n\big)\varphi(L_{-n},M_{n-2\mu}),\label{corII.1adn}\\
&&\!\!\!\!\!\!\!\!\!\!\!\!\!\!\!\!\!\!\!
(2n+m+2\mu) \varphi(L_{-m-2\mu},M_{m})\nonumber\\
&&\!\!\!\!\!\!\!\!\!\!\!\!\!\!\!\!\!\!\!
=(m+2\mu-n\lambda)\varphi(L_{-n-m-2\mu},M_{m+n})
-\big((n+m+2\mu)(1+\lambda)-n\big)\varphi(L_{n},M_{-n-2\mu}).\label{corII.1adm+n}
\end{eqnarray}

Then using the same arguments given in Lemma 2.2 of the reference
\cite{LS3}, we obtain the following results.

\ni If $\lambda=0$, then
\begin{eqnarray}\label{IISt3al1}
\varphi(L_{-m-2\mu},M_m) =\frac{(m+2\mu)(m+2\mu+1)}{2}
\varphi(L_{-1},M_{1-2\mu}),\ \ \forall\,\,m\in\Z.
\end{eqnarray}
If $\lambda=1$, then
\begin{eqnarray}\label{IISt3al2}
\varphi(L_{-m-2\mu},M_{m})=(m+2\mu-1)(m+2\mu)(m+2\mu+1)c'_1,\ \
\forall\,\,m\in\Z.
\end{eqnarray}
for some constant $c'_1\in\C$.

\ni If $\lambda=-1$, then
\begin{eqnarray*}
\varphi(L_{-m-2\mu},M_{m})=\left\{\begin{array}{ll}
(m+2\mu)\varphi(L_0,M_{-2\mu})&{\rm if}\ \,m\geq-2\mu,\vs{8pt}\\
(m+2\mu+2)\varphi(L_0,M_{-2\mu})&{\rm if}\ \,m<-2\mu.
\end{array}\right.
\end{eqnarray*}
by (\ref{corIIYYLM-1}), which can be rewritten as
\begin{eqnarray}\label{IISt3al4}
\varphi(L_{-m-2\mu},M_{m})=0,\ \ \forall\,\,m\in\Z.
\end{eqnarray}
If $\lambda\notin\{-1,\,0,\,1\}$, then
\begin{eqnarray}\label{IISt3al3}
\varphi(L_{-m-2\mu},M_{m})=0,\ \ \forall\,\,m\in\Z.
\end{eqnarray}

\ni{\bf Step 4.}\ \,The computation of
$\varphi(Y_{-2\mu-p},Y_p),\,\,\forall\,\,p\in\frac{1}{2}+\Z$.

Using the Jacobian identity on
$(L_{-m},\,Y_{-n-\mu},\,Y_{m+n-\mu})$, $\forall\,\,m,n\in\Z$, one has
\begin{eqnarray}
&&\big(m(1+\lambda)-2n\big)\varphi(Y_{-m-n-\mu},Y_{m+n-\mu})\nonumber\\
&&=2(m+2n)\varphi(L_{-m},M_{m-2\mu})-\big(m(3+\lambda)+2n\big)
\varphi(Y_{-n-\mu},Y_{n-\mu}).\label{corII.2lad-1}
\end{eqnarray}
In particular, taking $n=0$ in (\ref{corII.2lad-1}), one has
\begin{eqnarray}\label{corII.2lad-1asd}
m(1+\lambda)\varphi(Y_{-m-\mu},Y_{m-\mu})
=2m\varphi(L_{-m},M_{m-2\mu})-m(3+\lambda)
\varphi(Y_{-\mu},Y_{-\mu}).
\end{eqnarray}

\ni If $\lambda=-1$, then (\ref{corII.2lad-1}) becomes
\begin{eqnarray*}
\varphi(L_{-m},M_{m-2\mu})=\varphi(Y_{-\mu},Y_{-\mu}),\ \
\forall\,\,m\in\Z^*,
\end{eqnarray*}
which together with (\ref{IISt3al4}), forces
$\varphi(Y_{-\mu},Y_{-\mu})=\varphi(L_0,M_{-2\mu})=0$. Then
(\ref{IISt3al4}) can be rewritten as
\begin{eqnarray}\label{0801231}
\varphi(L_{-m-2\mu},M_{m})=0,\ \ \forall\,\,m\in\Z.
\end{eqnarray}
So in this case, (\ref{corII.2lad-1}) can be rewritten as
\begin{eqnarray*}
&&\varphi(Y_{-m-\mu},Y_{m-\mu})=m\varphi(Y_{-1-\mu},Y_{1-\mu}),\ \
\forall\,\,m\in\Z.
\end{eqnarray*}
Then (\ref{corII.2lad-1}) becomes\,
$m(m+n)\varphi(Y_{-1-\mu},Y_{1-\mu})=0$, which implies
$\varphi(Y_{-1-\mu},Y_{1-\mu})=0$, and further
\begin{eqnarray}\label{0801232}
&&\varphi(Y_{-m-\mu},Y_{m-\mu})=0,\ \ \forall\,\,m\in\Z.
\end{eqnarray}

\ni If $\lambda=0$, using (\ref{IISt3al1}) and
(\ref{corII.2lad-1asd}), we obtain
\begin{eqnarray}\label{IISt3al1ax}
&&\varphi(Y_{-m-\mu},Y_{m-\mu})=m(m+1)\varphi(L_{-1},M_{1-2\mu})
-3\varphi(Y_{-\mu},Y_{-\mu}),\ \ \forall\,\,m\in\Z^*.
\end{eqnarray}
which together with (\ref{corII.2lad-1}), gives\,\,
$6\varphi(Y_{-\mu},Y_{-\mu})+mn\varphi(L_{-1},M_{1-2\mu})=0,\,\,\forall\,\,n\in\Z,m\in\Z^*$
and further forces
$\varphi(Y_{-\mu},Y_{-\mu})=\varphi(L_{-1},M_{1-2\mu})=0$. Then
recalling (\ref{IISt3al1}) and (\ref{IISt3al1ax}), we obtain
\begin{eqnarray}\label{Ila1x}
\varphi(L_{-m-2\mu},M_m)=\varphi(Y_{-m-\mu},Y_{m-\mu})=0,\ \
\forall\,\,m,n\in\Z.
\end{eqnarray}

\ni If $\lambda=1$, using (\ref{IISt3al2}) and
(\ref{corII.2lad-1asd}), we obtain
\begin{eqnarray}
&&\varphi(Y_{-m-\mu},Y_{m-\mu})=m(m^2-1)c_1-2
\varphi(Y_{-\mu},Y_{-\mu}),\ \
\forall\,\,m\in\Z^*,\label{IISt3al2s2x}
\end{eqnarray}
for some constant $c_1\in\C$. Using (\ref{IISt3al2}),
(\ref{IISt3al2s2x}) and (\ref{corII.2lad-1}), we obtain
$\varphi(Y_{-\mu},Y_{-\mu})=0$. Then (\ref{IISt3al2s2x}) can be
rewritten as
\begin{eqnarray}
&&\varphi(Y_{-m-\mu},Y_{m-\mu})=m(m^2-1)c'_1,\ \
\forall\,\,m\in\Z^*.\label{IISt2xuy}
\end{eqnarray}

\ni If $\lambda\notin\{-1,\,0,\,1\}$, taking $n=0$ in using
(\ref{IISt3al3}) and (\ref{corII.2lad-1asd}), we obtain
\begin{eqnarray*}
&&\!\!\!\!\!\!\!\!\!m(1+\lambda)\varphi(Y_{-m-\mu},Y_{m-\mu})=-m(3+\lambda)
\varphi(Y_{-\mu},Y_{-\mu}),
\end{eqnarray*}
which gives
\begin{eqnarray}\label{IISt3al3lbs}
\varphi(Y_{-m-\mu},Y_{m-\mu})=-\frac{3+\lambda}{1+\lambda}
\varphi(Y_{-\mu},Y_{-\mu}),\ \ \forall\,\,m\in\Z^*.
\end{eqnarray}
Taking (\ref{IISt3al3lbs}) back to (\ref{corII.2lad-1}), one has
$\varphi(Y_{-\mu},Y_{-\mu})=0$, which gives
\begin{eqnarray}\label{IIsSt3al3lb}
\varphi(Y_{-m-\mu},Y_{m-\mu})=0,\ \ \forall\,\,m\in\Z.
\end{eqnarray}
Then this lemma follows.\QED

\ni Then the lemma proves Theorem \ref{mainth} (ii)--(v).

\begin{case}
\adddot\label{cs4} $\mu\in\Z$.
\end{case}
\begin{lemm}\adddot\label{lemm4}
(i)\ For the subcase $\lambda\ne-1$, one has $\varphi=0$.\\
(ii)\ For the subcase $\lambda=-1$, only
$\varphi(Y_{-p},Y_{p-2\mu})\,\,(\,\forall\,\,p\in\frac{1}{2}+\Z)$ is
not vanishing, given in (\ref{08012321}).
\end{lemm}
{\it Proof.} One has
\begin{eqnarray}
&&\varphi(L_0,Y_p)=0,\ \ \forall\,\,p\in\frac{1}{2}+\Z,\label{recIII2}\\
&&\varphi(L_0,M_n)=0\ \ {\rm if}\ \,n\neq-2\mu,\ \
\varphi(L_1,M_{-2\mu-1})=0\ \ {\rm if}\
\,\lambda\neq-1.\label{recIII3}
\end{eqnarray}

For any $m,n\in\Z,\,p,q\in\frac{1}{2}+\Z$, according to the
identities (\ref{corIIYM})--(\ref{corIILY}), (\ref{recIII2}) and
(\ref{recIII3}), one has
\begin{eqnarray}
&&(p+q+2\mu)\varphi(Y_p,Y_q)=(m+n+2\mu)\varphi(L_m,M_n)=0,\label{corIIILM}\\
&&\varphi(L_0,M_{-2\mu})
=(m+n+4\mu)\varphi(M_m,M_n)=\varphi(Y_p,M_m)=\varphi(L_m,Y_p)=0.\label{corIIIYM}
\end{eqnarray}
Similarly, one also can prove (\ref{corIIMMla}) holds in this case.
Then from (\ref{recIII2})--(\ref{corIIIYM}), the left components we
have to present in this case are listed in the following (\,where
$m$ is an arbitrary integer):
\begin{eqnarray*}
&&\varphi(Y_{-2\mu-m},Y_m)\,\,{\rm and}\,\,\varphi(L_{-2\mu-m},M_m),
\end{eqnarray*}
which will be taken into account together in the following.

As in Case \ref{cs2}, the results (\ref{IISt3al1})--(\ref{IISt3al3})
still hold in this case. Hence we shall cite them directly here
using the same notations. Using the Jacobian identity on
$(L_{-m},\,Y_{-p},\,Y_{m+p-2\mu})$,
$\forall\,\,m\in\Z,\,p\in\frac{1}{2}+\Z$, one has
\begin{eqnarray}\label{corIII.2lad-1}
&&\big(2\mu-2p+m(1+\lambda)\big)\varphi(Y_{-m-p},Y_{m+p-2\mu})\nonumber\\
&&=2(m+2p-2\mu)\varphi(L_{-m},M_{m-2\mu})
-\big(2p-2\mu+m(3+\lambda)\big)\varphi(Y_{-p},Y_{p-2\mu}).
\end{eqnarray}

\ni If $\lambda=-1$, then (\ref{corIII.2lad-1}) becomes
\big(\,recalling(\ref{IISt3al4})\big)
\begin{eqnarray*}
&&(\mu-p)\varphi(Y_{-m-p},Y_{m+p-2\mu})+(p-\mu+m)\varphi(Y_{-p},Y_{p-2\mu})=0,
\end{eqnarray*}
from which we can deduce
\begin{eqnarray}\label{08012321}
\varphi(Y_{-p},Y_{p-2\mu})=(p-\mu)c,\ \
\forall\,\,p\in\frac{1}{2}+\Z.
\end{eqnarray}
for some constant $c\in\C$.

\ni If $\lambda=0$, using (\ref{IISt3al1}) and
(\ref{corIII.2lad-1}), we obtain (\,taking $m=-1$)
\begin{eqnarray}\label{0801241}
\big(2(p+1)-(2\mu+1)\big)\varphi(Y_{1-p},Y_{p-1-2\mu})
=\big(2(p-1)-(2\mu+1)\big)\varphi(Y_{-p},Y_{p-2\mu}).
\end{eqnarray}
Using induction on $n=p-\frac{1}{2}$ where $p$ is that determined in
(\ref{0801241}), one has
\begin{eqnarray}\label{0801242}
&&\varphi(Y_{-p},Y_{p-2\mu})=\big(2(p+1)-(2\mu+1)\big)\big(2p-(2\mu+1)\big)c',
\end{eqnarray}
for some constant $c'\in\C$. Taking $\varphi(L_{-m-2\mu},M_m)$ and
$\varphi(Y_{-p},Y_{p-2\mu})$ respectively given in (\ref{IISt3al1})
and (\ref{0801242}) back to the identity (\ref{corIII.2lad-1}), we
obtain
\begin{eqnarray*}
m(m+1)\big(4(m-1)\varphi(L_{-1},M_{1-2\mu})+(2p+2\mu-m)c'\big)=0,\ \
\forall\,\,m\in\Z,\,p\in\frac{1}{2}+\Z,
\end{eqnarray*}
which forces $\varphi(L_{-1},M_{1-2\mu})=c'=0$. Then
(\ref{IISt3al1}) and (\ref{0801242}) becomes
\begin{eqnarray}\label{0801243}
&&\varphi(L_{-m-2\mu},M_m)=\varphi(Y_{-p},Y_{p-2\mu})=0,\ \
\forall\,\,m\in\Z,\,p\in\frac{1}{2}+\Z.
\end{eqnarray}

\ni If $\lambda=1$, using (\ref{IISt3al2}) and
(\ref{corIII.2lad-1}), we obtain (\,taking $m=-1$)
\begin{eqnarray}\label{10801244}
&&(p+1-\mu)\varphi(Y_{1-p},Y_{p-1-2\mu})
=(p-2-\mu)\varphi(Y_{-p},Y_{p-2\mu}).
\end{eqnarray}
Using induction on $n=p-\frac{1}{2}$ where $p$ is that given in
(\ref{10801244}), one has
\begin{eqnarray}\label{10801245}
&&\varphi(Y_{-p},Y_{p-2\mu})=(p-1-\mu)(p-\mu)(p+1-\mu)c'',\ \
\forall\,\,p\in\frac{1}{2}+\Z.
\end{eqnarray}
Taking $\varphi(L_{-m-2\mu},M_{m})$ and $\varphi(Y_{-p},Y_{p-2\mu})$
respectively given in (\ref{IISt3al2}) and (\ref{10801245}) back to
(\ref{corIII.2lad-1}), one can deduce
\begin{eqnarray*}
&&c'_1m(m^2-1)(m+2p-2\mu)-c''m\big(m^2+p^2+m(p-\mu-1)-2p(1+\mu)+\mu(2+\mu)\big),
\end{eqnarray*}
for any  $m\in\Z,\,p\in\frac{1}{2}+\Z$, which forces $c'_1=c''=0$.
Then
\begin{eqnarray}\label{10801246}
&&\varphi(L_{-m-2\mu},M_{m})=\varphi(Y_{-p},Y_{p-2\mu})=0,\ \
\forall\,\,m\in\Z,\,p\in\frac{1}{2}+\Z.
\end{eqnarray}

\ni If $\lambda\notin\{-1,\,0,\,1\}$, using (\ref{IISt3al3}), one
can rewrite (\ref{corIII.2lad-1}) as (\,taking $m=1$)
\begin{eqnarray*}
&&(2p-2\mu-1-\lambda)\varphi(Y_{-1-p},Y_{p+1-2\mu})
=(2p-2\mu+3+\lambda)\varphi(Y_{-p},Y_{p-2\mu}).
\end{eqnarray*}
Similarly, one also can deduce (\ref{IIsSt3al3lb}) holds in this
case. Then this lemma follows.\QED

\ni Then the lemma proves Theorem \ref{mainth} (vi).

\end{document}